\documentclass[11pt, a4paper]{article}

\usepackage{authblk}

\usepackage[a4paper]{geometry}
\usepackage{lineno}

\usepackage{amsmath}
\usepackage{amssymb}
\usepackage{amsthm}
\usepackage{amsfonts}
\usepackage{dsfont}
\usepackage{mathtools}

\usepackage{booktabs}

\usepackage{algorithm}
\usepackage{algpseudocode}

\usepackage{graphicx}
\usepackage{caption}

\usepackage{parskip}

\usepackage{xcolor}

\usepackage{lmodern}

\usepackage{tikz}
\usetikzlibrary{arrows.meta, positioning}

\newcommand{\R}{\mathbb{R}}
\newcommand{\E}{\mathbb{E}}
\newcommand{\la}{\langle}
\newcommand{\ra}{\rangle}

\begin{document}

\title{Variance Reduction in the Fokker--Planck Particle Method for Rarefied Gases using Quasi-Random Numbers}
\author[1]{Lukas Netterdon \thanks{netterdon@acom.rwth-aachen.de}}
\author[2, 3]{Veronica Montanaro\thanks{veronica.montanaro@empa.ch}}
\author[1]{Manuel Torrilhon \thanks{mt@acom.rwth-aachen.de}}
\author[3]{Hossein Gorji \thanks{mohammadhossein.gorji@empa.ch}}

\affil[1]{Chair of Applied and Computational Mathematics, RWTH Aachen University, D-52062 Aachen, Germany}
\affil[2]{Chair of Scientific Computing and Uncertainty Quantification, \'Ecole Polytechnique F\'ed\'erale de Lausanne, CH-1015 Lausanne, Switzerland}
\affil[3]{Laboratory for Computational Engineering, Swiss Federal Institute of Material Sciences, CH-8600 D\"ubendorf, Switzerland}

\maketitle

\begin{abstract}

The Fokker--Planck (FP) particle method accelerates rarefied-gas simulations by replacing the binary collisions of the commonly used Direct Simulation Monte Carlo (DSMC) method with a drift--diffusion process. Like all particle methods, the FP method is inherently stochastic, which leads to statistical fluctuations in macroscopic quantities and necessitates large particle numbers for accurate results. In this work, we investigate the use of quasi-random numbers, which sample distributions more evenly and thereby reduce the variance. To preserve the low-discrepancy structure across time steps, we employ the Array Randomized Quasi-Monte Carlo (Array-RQMC) technique. We combine the FP method with Array-RQMC and compare it in homogeneous and inhomogeneous problems with other commonly used variance-reduction techniques. The proposed FP--Array-RQMC approach achieves improved convergence rates compared with pseudo-random sampling and yields smaller estimator errors for sufficiently large particle numbers.
\end{abstract}

\section{Introduction}
For rarefied gases, the mean free path $\lambda$ is large compared to a characteristic length $L$, i.e., the Knudsen number $Kn=\lambda/L$ is not negligible. Such gases can deviate significantly from thermal equilibrium and are therefore not accurately described by continuum models such as the Navier--Stokes--Fourier equations \cite{struchtrup_macroscopic_2005}. Instead, a statistical description is required in which molecular interactions are taken into account. In most cases, these interactions can be modeled appropriately as binary collisions, leading to the Boltzmann equation \cite{cercignani_boltzmann_1988}.

Since the statistical description leads to a high-dimensional phase space, direct simulations of the Boltzmann equation are computationally infeasible. In practice, the particle method Direct Simulation Monte Carlo (DSMC) is therefore used, which advances an ensemble of particles by alternating free-flight and binary-collision steps \cite{bird_molecular_1994}. However, in low-Knudsen-number regimes, where the collision frequency is high, the explicit treatment of binary collisions becomes prohibitive.

The Fokker--Planck (FP) particle method was introduced to accelerate simulations by replacing explicit binary collisions with a drift--diffusion process. In this formulation, particles evolve along stochastic paths that are coupled only through velocity moments, which decouples the computational complexity from the collision frequency. The drift and diffusion coefficients are chosen such that the Boltzmann equation is approximated in a manner that preserves its key physical properties \cite{jenny_solution_2010, gorji_entropic_2021, montanaro_fisher_2025}.

Both DSMC and the FP method are inherently stochastic, so macroscopic quantities estimated from finite particle ensembles exhibit statistical fluctuations. A natural goal is to reduce these fluctuations without increasing the number of particles. Since the FP method is formulated in terms of stochastic differential equations, it is well suited for variance-reduction techniques. In this work, we replace the pseudo-random numbers used in the velocity update by quasi-random numbers in order to reduce the estimator variance.

Quasi-random numbers sample distributions more evenly than pseudo-random numbers and therefore improve the accuracy of moment estimates \cite{niederreiter_random_1992}. When applied to stochastic processes, however, care is required to preserve their low-discrepancy structure across successive steps. The Array Randomized Quasi-Monte Carlo (Array-RQMC) method addresses this challenge by reconstructing the low-discrepancy structure at each step by sorting the particle states \cite{lecuyer_sorting_2018}.

In this work, we integrate the Array-RQMC technique into the FP particle method to exploit the advantages of quasi-random numbers in a spatially discretized setting with particle transport across cells and boundary interactions. We consider both homogeneous and inhomogeneous test cases and compare the proposed approach with pseudo-random sampling and other commonly used variance-reduction techniques. The numerical experiments show that the method yields improved convergence rates and smaller estimator errors for sufficiently large particle numbers.

The paper is structured as follows. Section~2 reviews the FP equation and the corresponding particle algorithm. Section~3 introduces quasi-random numbers and the Array-RQMC method. Section~4 presents their integration into the FP method. Section~5 reports the numerical results, and Section~6 summarizes the main findings.

\section{Fokker-Planck Method}
Gases are described statistically by a molecular velocity distribution
\begin{equation}
f(c,x,t)=\rho(x,t)\hat f(c|x,t)
\end{equation}
where $\rho$ is the mass density and $\hat f$ is the conditional probability density function of the velocities $c\in\R^3$ at position $x\in\R^3$ at time $t\in\R^+$.

Macroscopic quantities are obtained as velocity moments of $f$, defined by
\begin{equation} \label{eq:moments}
   \la g(c)|f \ra = \int_{\R^3}g(c)f(c,x,t)dc .
\end{equation}
Important moments are the mass density $\rho$, the momentum density $\rho\bar c_i$ (with bulk velocity $\bar c_i$), the translational energy density $\epsilon$, the stress tensor $\sigma_{ij}$, and the heat flux $q_i$. These are obtained by choosing
\begin{equation}
g(c)\in\{1, c_i, \frac12  u_ku_k, u_{\la i}u_{j\ra}, \frac12 u_i u_k u_k\},
\end{equation}
respectively. Here, $u_i = c_i - \bar c_i$ are the peculiar velocities.  The angular bracket in $u_{\la i}u_{j\ra}$ denotes the symmetric, trace-free part of the tensor. Throughout, the Einstein summation convention is used for repeated indices.    

\subsection{Kinetic Equations}
The evolution of the distribution function is governed by a kinetic equation of the form
\begin{equation} 
  \frac{\partial}{\partial t}(f) + \frac{\partial}{\partial x_i} (c_i f) = S(f), 
\end{equation} 
where the left-hand side describes free flight and $S(f)$ denotes the collision operator.

The collision operator must satisfy three key properties:
\begin{enumerate}
\item The distribution relaxes toward the correct equilibrium, given by the Maxwell distribution, 
\begin{equation}
f^{Maxwell} = \rho \Lambda \exp\left(-\frac{(c_i-\bar c_i)^2}{2kT/m}\right),
\end{equation}  
where $\Lambda$ is a normalization factor, $k$ is the Boltzmann constant, $m$ the molecular mass and $T = \frac23 \frac{m}{k} \epsilon$ the temperature.
\item The H-theorem is fulfilled, ensuring that the physical entropy is non-decreasing.
\item The relaxation rates of relevant moments are correct:
\begin{equation}  
P_{\la g(c)|f \ra} = \int_{\R^3} g(c) S(f) dc,
\end{equation}  
where $g(c)=\{1, c_i, \frac12  u_ku_k, u_{\la i}u_{j\ra}, \frac12 u_i u_k u_k\}$ correspond to the relaxation rates for mass density, momentum density, energy density, stress tensor, and heat flux. 
\end{enumerate}

\subsubsection{Boltzmann Equation}
When collisions are modeled as binary interactions between particles, the collision operator is given by the Boltzmann equation \cite{struchtrup_macroscopic_2005}:
\begin{equation}
S^{\text{Boltz}}(f) = \int_{\R^3}\int_0^{2\pi}\int_0^{\pi/2} (f'f'_1-ff_1)c_r\sigma \sin(\theta) d\theta d\phi dc_1. 
\end{equation}
Here, $f=f(c,x,t)$ and $f_1=f(c_1,x,t)$ denote the pre-collision distributions of the collision pair $c$, $c_1$, while $f'$ and $f_1'$ denote the corresponding post-collision distributions. $c_r = |c -c_1|$ is the relative velocity, $\sigma$ is the differential cross-section, $\theta$ the scattering angle, and $\phi$ the orientation of the collision plane. 

The Boltzmann equation drives the distribution function toward the Maxwell distribution and satisfies the H-theorem. Numerous numerical studies demonstrate that the relaxation rates of relevant moments are captured accurately across a wide range of scenarios \cite{cercignani_boltzmann_1988}. For this reason, the Boltzmann equation serves as the reference model against which other kinetic models are compared.

However, due to the high dimensionality of the distribution function, algorithms requiring a discretization of the (x-c-t)-space are prohibitively computationally expensive. In practice, the Boltzmann equation is typically solved using a particle method called Direct Simulation Monte Carlo \cite{bird_molecular_1994}. Since DSMC computes binary collisions explicitly, it becomes computationally expensive, especially in the low Knudsen-number regime.

\subsubsection{Fokker-Planck Equation}
To improve computational efficiency, alternative kinetic models have been considered. One approach describes collisions as a drift-diffusion process, leading to the Fokker-Planck equation \cite{kirkwood_statistical_1946, pawula_approximation_1967}:
\begin{equation}
 S^{\text{FP}}(f) = - \frac{\partial}{\partial c_i}(A_i(f)f) + \frac{\partial^2}{\partial c_i \partial c_i} (D(f) f), 
\end{equation}
where $A_i$ is the drift and $D$ the diffusion coefficient. Both depend on velocity moments and must be chosen to satisfy the three properties above. In particular, the relaxation rates of the FP equation must match those of the Boltzmann equation for relevant moments: 
\begin{equation}
P_{\la g(c)|f \ra}^{\text{FP}} = \int_{\R^3} g(c) S^{\text{FP}}(f) dc = \int_{\R^3} g(c) S^{\text{Boltz}}(f) dc = P_{\la g(c)|f \ra}^{\text{Boltz}}.
\end{equation}
In this sense, the $FP$ equation may be interpreted as an approximation of the Boltzmann equation.

The simplest viable model uses a linear drift and a scalar diffusion \cite{jenny_solution_2010}.
\begin{align} \label{eq:drift_diffusion}
&A_i(f) = \frac{1}{\tau}(c_i-\bar c_i), \\
&D(f) = \frac{2\epsilon}{3\tau},
\end{align}
with free parameter $\tau$. This model satisfies conservation laws, the H-theorem, and relaxes to the correct equilibrium distribution. The parameter $\tau$ is chosen to ensure the correct relaxation rate for the pressure tensor. However, because there are no further free parameters, the relaxation rate of the heat flux cannot be adjusted, leading to an incorrect Prandtl number. To address this limitation, more advanced FP models have been developed \cite{gorji_fokkerplanck_2011, mathiaud_fokkerplanck_2017, gorji_entropic_2021, montanaro_fisher_2025}.

Since this paper focuses on variance reduction, the linear FP model will be sufficient for our purposes. 

\subsection{Fokker-Planck Algorithm}
Although the FP collision operator is simpler than the Boltzmann operator, the high dimensionality of the distribution function still makes direct simulation infeasible. Instead, the gas is modeled using an ensemble of representative particles, whose trajectories follow stochastic paths coupled only through velocity moments. 

The particle velocities $M_i(t)$ evolve according to the Langevin equation:
\begin{equation}\label{eq:langevin}
dM_i(t) = A_i[f](X, M)dt + \sqrt{2D[f](X, M)} dW_i,  
\end{equation}
where $dW_i$ is a Wiener process with zero mean and variance $dt$. The corresponding particle positions $X_i(t)$ evolve as
\begin{equation}\label{eq:free_flight}
dX_i(t) = M_i dt.
\end{equation}
The joint density of the random variables $M$, $X$, evolving under this Langevin dynamics, coincides with the solution of the FP equation \cite{risken_fokker-planck_1996}. 

Time integration of this system requires discretization in both space and time, since drift and diffusion depend on macroscopic moments. In each time step, the moments are computed, drift and diffusion coefficients evaluated, and the particle velocities and positions updated. The algorithm can be interpreted as a variant of DSMC in which explicit binary collisions are replaced by stochastic velocity updates. The key advantage lies in computational efficiency, because the velocity update scales linearly with the number of particles and does not scale with collision frequency. This makes the FP method especially well-suited for low-Knudsen-number regime \cite{gorji_efficient_2014, gorji_fokkerplanckdsmc_2015}.

Both DSMC and the FP method are stochastic in nature, i.e., their updates rely on random variables that introduce statistical fluctuations into macroscopic quantities. These fluctuations can be reduced by increasing the number of particles, but this comes at a higher computational cost. Since the FP method evolves particles along continuous stochastic paths rather than discrete collisions, it admits the use of variance reduction techniques developed for stochastic differential equations. One approach is to replace pseudo-random numbers by quasi-random numbers, which are deterministic sequences designed to sample the distributions more evenly. As a result, estimates converge faster to their expected values, with reduced fluctuations. \\
The next chapter introduces the concept of quasi-random numbers in more detail and explains how they can be integrated into stochastic particle simulations through the Array Randomized Quasi-Monte Carlo method.  
  
\section{Array Randomized Quasi-Monte Carlo Method}
Let $Y\in\R^d$ be a random variable with probability density function $p_Y(y)$. We are often interested in the expected value of a measurable function $g: \R^d \to \R$,
\begin{equation}\label{eq:mu}
    \mu = \mathbb{E}\!\left[g(Y)\right] = \int_{\R^d} g(y)p_Y(y) dy.
\end{equation}
In most cases, the expectation $\mu$ cannot be evaluated analytically because the density $p_Y(y)$ is only implicitly available or the integrand $g(y)p_Y(y)$ is too complex to integrate in closed form. Instead, the Monte Carlo (MC) method is used to approximate $\mu$ through $N$ independent realizations $\{Y^{(i)}\}_{i=1}^N$ of $Y$, yielding the estimator
\begin{equation}\label{eq:hat_mu}
\hat\mu_{N} = \frac{1}{N}\sum_{i=1}^N g\!\left(Y^{(i)}\right).
\end{equation}
According to the law of large numbers, $\hat\mu_{N} \to \mu$ as $N \to \infty$, and by the central limit theorem, the root-mean-square error decreases as \cite{gentle_random_2003}
\begin{equation}
    \mathrm{RMSE}(\hat\mu_{N}) = O(N^{-1/2}).
\end{equation}
This dimension-independent convergence rate is the main strength and the main limitation of MC methods, since it is robust but relatively slow \cite{robert_monte_2004}. This motivates the use of quasi-MC (QMC) techniques, which aim to improve convergence by replacing independent random numbers with more evenly distributed point sets.

\subsection{Quasi-Random Numbers}
Samples from any probability distribution can be generated from uniform random variables over $[0,1]^d$ by applying an appropriate transformation, such as the inverse transform method in one dimension or more general sampling schemes in higher dimensions \cite{devroye_non-uniform_1986}. Consequently, the study of sampling strategies reduces to the placement of points in the unit cube. Point sets which are more evenly distributed in $[0,1]^d$ yield smaller integration errors and therefore reduce the variance of estimators.

The notion of evenness is formalized by discrepancy, which measures the deviation of a point set from perfect uniformity \cite{niederreiter_random_1992, dick_digital_2010}. Let $P_N = \{u^{(1)}, ..., u^{(N)}\} \subset [0,1]^d$ be a finite point set in the $d$-dimensional unit cube. For any axis-aligned box anchored at the origin,
\begin{equation}
B(x) = [0,x) = \prod_{j=1}^d [0, x_j], \quad x \in [0,1]^d,
\end{equation}
the proportion of points of $P_N$ contained in $B(x)$ is
\begin{equation}
\frac1N \sum_{i=1}^N \mathds{1}_{B(x)}(u^{(i)}).
\end{equation}
Here $\mathds{1}_{B(x)}$ is an indicator function. The star discrepancy of $P_N$ is defined as
\begin{equation}
D^*(P_N) = \sup_{x\in[0,1]^d}\left|\frac1N \sum_{i=1}^N \mathds{1}_{B(x)}(u^{(i)}) - \lambda(B(x))\right|,
\end{equation}
where $\lambda(B(x)) = \prod_{j=1}^d x_j$ is the Lebesgue measure of the box. \\
The star discrepancy bounds the integration error when approximating an integral by the average of function evaluations at the points $P_N$. The Koksma-Hlawka inequality states that if $f:[0,1]^d \to \R$ has a bounded variation in the sense of Hardy-Krause, then   
\begin{equation}  
\left|\frac1N \sum_{i=1}^N f(u^{(i)}) - \int_{[0,1]^d}f(u)du\right| \leq V(f)D^*(P_N), 
\end{equation}  
where $V(f)$ denotes the Hardy-Krause variation of $f$.

Low-discrepancy point sets can be constructed explicitly. A popular choice is given by Sobol' sequences, since they are efficient to generate, perform well in higher dimensions, and are implemented professionally in many software libraries \cite{sobol_distribution_1967, joe_constructing_2008}. An important property of Sobol' sequences is that they can be generated sequentially and therefore every consecutive block of points retains good uniformity properties.

Sobol' sequences have a provable bound on their evenness. Their star discrepancy satisfies
\begin{equation}
D^*(P_N^{(\text{Sobol})}) = O\left(\frac{(\log N)^d}{N}\right),
\end{equation}

which implies through the Koksma--Hlawka inequality that the integration error is bounded by the same order for functions of bounded variation. In contrast, the star discrepancy of independent random point sets is $D^*(P_N^{(\text{Random})}) = O(N^{-1/2})$. Thus, Sobol' sequences achieve a strictly better asymptotic error bound than independent random sampling.

Since low-discrepancy sequences are constructed deterministically, it is not directly possible to calculate the variance of an estimator, which complicates the analysis of convergence rates. To overcome this limitation, randomized QMC (RQMC) methods apply a randomization to the low-discrepancy set, in such a way that the points remain uniformly distributed and still preserve the low-discrepancy structure. This yields estimators with well-defined variance, so that the convergence rates of RQMC and MC can be compared, and in general RQMC achieves faster convergence than MC \cite{owen_monte_1997, lecuyer_randomized_2018}.

\begin{figure}[h!]
\centering
\includegraphics[width=1\linewidth]{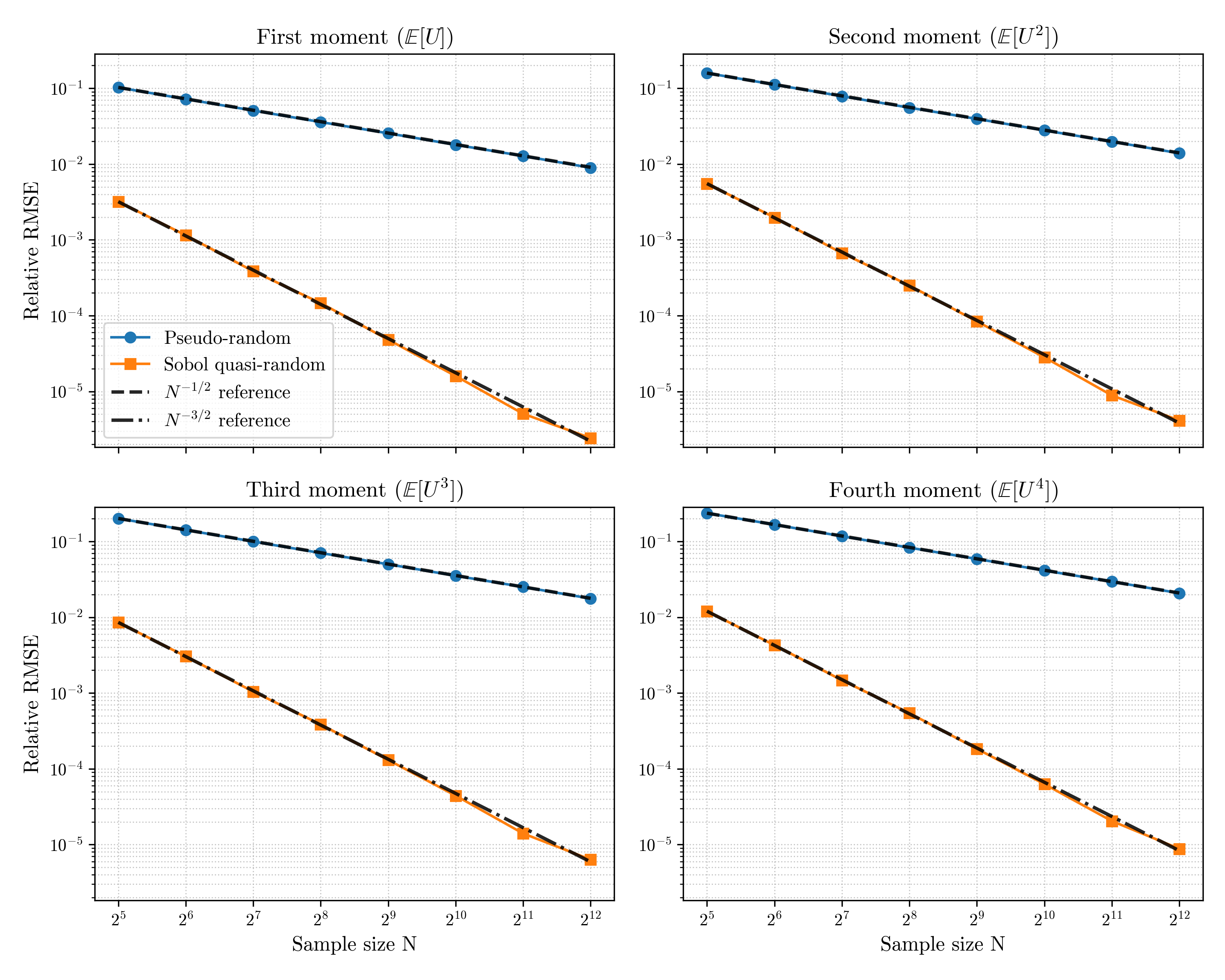}
\caption{Convergence of MC and QMC estimators for the first four moments of a uniform random variable. The plot shows the relative RMSE of the estimated moments versus the sample size $N$ on a log-log scale. Reference slopes proportional to $N^{-1/2}$ and $N^{-3/2}$ are added for comparison. QMC consistently outperforms MC.} 
\label{fig:uniform_moment_convergence}
\end{figure}

To illustrate the impact of sampling more evenly, Figure \ref{fig:uniform_moment_convergence} compares the convergence of MC and RQMC estimators for a simple one-dimensional test case. Uniform samples $U_i \in (0,1)$ are used to estimate the first four moments of the uniform distribution, 
\begin{equation}
\mu_k = \E[U^k] = \frac{1}{k+1}, \quad k=1,...,4
\end{equation}
The estimators
\begin{equation}
\hat \mu_{k,N} = \frac{1}{N} \sum_{i=1}^N U_i^k
\end{equation}
are computed from both pseudo-random and randomized Sobol' sequences. The relative RMSE is calculated over many repetitions. \\
In the MC case, the error decreases as expected with $O(N^{-1/2})$. In the QMC case, the error decreases with $O(N^{-3/2})$, in agreement with theoretical results for randomized low-discrepancy sequences applied to smooth integrands. \\

This example illustrates the improved accuracy and faster convergence that motivates the use of quasi-random numbers in stochastic particle methods. However, when extending these ideas to stochastic processes, additional care is required to preserve the low-discrepancy structure across time steps. 

\subsection{Array-RQMC}
Many stochastic problems of interest, such as particles whose trajectories are governed by the Langevin equation, can be represented as discrete-time Markov chains.

Let $Y_t \in \mathbb{R}^d$ denote the state of the system at time $t$. Its evolution is specified by a measurable mapping $\phi:\mathbb{R}^d \times (0,1)^d \to \mathbb{R}^d$ according to
\begin{equation}
    Y_{t+dt} = \phi(Y_t, U_t),
\end{equation}
where $U_t$ are independent and identically distributed (i.i.d.) uniform random variables over the unit hypercube $(0,1)^d$. In general, the distribution $Y_t$ is not available in closed form. Therefore, the MC method is used to investigate how an ensemble of $N$ independent realizations $\{Y_0\}_{i=1}^N$ progresses by updating the state successively with $N$ uniform realizations $\{U_t\}_{i=1}^N$, 
\begin{equation}\label{eq:numerical_transition_rule}
    Y^{(i)}_{t+dt} = \phi(Y^{(i)}_t, U^{(i)}_t), \quad i = 1, ..., N.
\end{equation}
As in the time-independent case, the desired expectation $\E[g(Y_t)]$ is estimated by averaging over the ensemble $\{Y_t\}_{i=1}^N$ (Eq.~\ref{eq:hat_mu}). To improve the accuracy of this estimation, we want to deploy quasi-random numbers. 

A naive way to apply RQMC to Markov chains is to replace the uniform samples $\{U_t\}_{i=1}^N$ in the transition rule (Eq.~\ref{eq:numerical_transition_rule}) with a low-discrepancy point set, meaning that, at each step, the current state $\{Y_t\}_{i=1}^N$ is paired with $\{U_t\}_{i=1}^N$. If $\{U_t\}_{i=1}^N$ is generated sequentially by a Sobol' sequence (no randomization), an artificial correlation between steps is introduced, distorting the dynamic of the chain. If, on the other hand, the set $\{U_t\}_{i=1}^N$ is randomized, the low-discrepancy structure between steps is lost, deteriorating the advantages of RQMC over MC.

 The Array-RQMC method was introduced to overcome this problem, by coupling the state and the uniform samples in an advantageous manner \cite{lecuyer_randomized_2005, lecuyer_sorting_2018}. To build intuition, we first consider the one-dimensional case. Assume that the number of realizations $N$ is large, so that the ensemble $\{Y_t\}_{i=1}^N$ can be depicted as a fine histogram. Realizations falling into the same narrow bin $[a-\epsilon, a+\epsilon]$ are nearly identical. Consequently, the updated ensemble of this bin $\{Y_{\text{bin},t+1}\}_{i=1}^N$ depends almost exclusively on the quality of the uniform realizations $\{U_{\text{bin},t}\}_{i=1}^N$, i.e., how evenly they are distributed in $(0,1)$. Quasi-random numbers are specifically designed to possess this property. The same procedure can be applied to all bins, to improve the accuracy of estimators for the entire ensemble.

The histogram perspective also clarifies why the naive application of RQMC to Markov chains fails to preserve the low-discrepancy structure between time steps: The uniform realizations that are coupled with the state in each bin, are random subsets from the entire low-discrepancy point set, and therefore lose the desired low-discrepancy properties.  

Let us now return to the Array-RQMC method. Constructing explicit histograms is cumbersome and requires choosing a bin width that may depend on the distribution. The Array-RQMC method captures the same idea without the need for a discretization. Instead of binning, the states $\{Y_t\}_{i=1}^N$ are sorted and then matched with a quasi-random point set. The sorting ensures that realizations with similar values are next to each other in the set. This is comparable to similar realizations being in the same bin. Since contiguous subsequences of low-discrepancy sequences approximately retain the low-discrepancy structure, state realizations with similar values will be updated with evenly distributed uniform realizations.   

To extend this procedure to higher dimensions, it is necessary to maintain the key property that in the ordered set $\{Y_t\}_{i=1}^N$, realizations that are next to each other, are also close in the $d$-dimensional state space. This is achieved by sorting the realizations along a space-filling curve. 
  
The benefits of the Array-RQMC method have been demonstrated in several applications, including option pricing \cite{abdellah_array-rqmc_2019} and stochastic chemical reaction networks \cite{puchhammer_variance_2021}. To the best of our knowledge, the method has not yet been applied in the domain of statistical mechanics, nor to spatially discretized multi-cell particle systems. In the next chapter, we describe how Array-RQMC is integrated into the Fokker--Planck particle method.

\section{Integrating Array-RQMC into the Fokker-Planck Method}
The FP method simulates rarefied gases by following the trajectories of an ensemble of representative particles $\{X^{(i)}, M^{(i)}\}$, which serve as a discrete approximation of the distribution function $f(x,c,t)$. The particle velocities evolve according to the Langevin equation (Eq.~\ref{eq:langevin}) and the positions by free flight (Eq.~\ref{eq:free_flight}). Because the drift and diffusion coefficients depend on macroscopic quantities, the method requires discretization in both space and time. The algorithm then proceeds as follows: 

\begin{algorithm}[H]
\caption{Fokker-Planck Method}
\begin{algorithmic}[1]
  \State \textbf{Input:} 
  initial particle ensemble $\{X^{(i)}, M^{(i)}\}$, 
  spatial grid $\{\Omega_j\}_{j=1}^{N_\text{cells}}$, 
  time step $\Delta t$, 
  total time $T$
  \For{$t_\text{idx} = 1, \dots, T/\Delta t$} 
    \For{each cell $\Omega_j$} 
      \State Compute macroscopic quantities
      \State Evaluate drift and diffusion coefficients
      \State Update particle velocities and positions
      \State Apply boundary conditions
    \EndFor
    \State Reassign particles that moved to neighboring cells
    \State Output macroscopic quantities as required
  \EndFor
\end{algorithmic}
\end{algorithm}
In the following subsections the steps in the algorithm are explained in more detail, beginning with the computation of macroscopic quantities. 

\subsection{Macroscopic Quantities}
Macroscopic quantities are obtained as moments of the distribution function (Eq.~\ref{eq:moments}). 
In the particle algorithm, these moments are estimated from the ensemble of representative particles. 
Let $\{X^{(i)}, M^{(i)}\}_{i=1}^{N_j}$ denote the particles located in cell $\Omega_j$ at a given time. 
The estimator of a moment of a function $g(M)$ is
\begin{equation}
\hat{\mu} = \frac{1}{N_j} \sum_{i=1}^{N_j} g(M^{(i)}),
\end{equation}
which serves as the discrete analog of $\langle g(c)\,|\,f\rangle$.
Of particular interest are the mean velocity, the translational energy, and the pressure tensor:
\begin{align}
\bar{c}_{k} &= \frac{1}{N_j} \sum_{i=1}^{N_j} M_{k}^{(i)}, \\
{\epsilon} &= \frac{1}{2N_j} \sum_{i=1}^{N_j} (M_{l}^{(i)}-\bar{c}_{l})(M_{l}^{(i)}-\bar{c}_{l}), \\
{\sigma}_{kl} &= \frac{1}{N_j} \sum_{i=1}^{N_j} (M_{k}^{(i)}-\bar{c}_{k})(M_{l}^{(i)}-\bar{c}_{l}) - \frac{1}{3}{\epsilon}\,\delta_{kl}.
\end{align}
Other macroscopic quantities, such as the heat flux, can be obtained analogously.

The drift and diffusion coefficients required in the Langevin equation are obtained from these macroscopic quantities (Eq.~\ref{eq:drift_diffusion}). Since they depend directly on ensemble estimates, statistical fluctuations can lead to inaccuracies in the computed coefficients. The use of quasi-random numbers in the proposed method is intended to reduce the variance of these estimators and thereby limit the propagation of noise.   

\subsection{Particle Update}
The evolution of the particle ensemble is governed by the Langevin equations (Eq.~\ref{eq:langevin}) together with the free flight motion (Eq.~\ref{eq:free_flight}). The discretization is performed in both space and time: particles in each cell evolve under locally constant drift and diffusion coefficients for a discrete time step $\Delta t$.

\subsubsection{Time Discretization}
For one velocity component, the linear FP model yields
\begin{equation}
dM
= -\frac{1}{\tau}\,(M-\bar c)\,dt
+ \sqrt{\frac{4\epsilon}{3\tau}}\,dW,
\end{equation}
A standard way of discretizing a stochastic differential equation is the Euler-Maruyama scheme \cite{kloeden_numerical_1992},
\begin{equation}
M_{t+\Delta t}
= M_t - \frac{\Delta t}{\tau}(M_t-\bar c)
+ \sqrt{\frac{4\epsilon}{3\tau}\,\Delta t}\,\xi_t,
\qquad \xi_t\sim\mathcal{N}(0,1).
\end{equation}
Since this method is only first-order accurate in time, relatively small time steps are required to obtain reliable results.
Instead, by recognizing the linear Langevin equation as an Ornstein-Uhlenbeck (OU) process, the known analytical solution can be used to remain accurate for any time step \cite{risken_fokker-planck_1996},
\begin{equation}
M_{t+\Delta t}
= \bar c + (M_t-\bar c)\,e^{-\Delta t/\tau}
+ \sqrt{\frac{2\epsilon}{3}\left(1-e^{-2\Delta t/\tau}\right)}\,\xi_t,
\qquad \xi_t\sim\mathcal{N}(0,1).
\end{equation}
Although the OU update eliminates the time-discretization error, the evolution of the particle ensemble remains stochastic through the Gaussian random variables $\xi_t$. These random samples determine the accuracy of the estimated macroscopic quantities and are therefore the target for variance reduction using quasi-random numbers.

\subsubsection{Fokker-Planck Method with Array-RQMC}
Array-RQMC requires a one-dimensional ordering of particle states such that nearby entries in the ordering correspond to similar states in velocity space. In the FP case, the relevant state variable is the three-dimensional velocity vector 
$M_t^{(i)} = (M_{t,x}^{(i)}, M_{t,y}^{(i)}, M_{t,z}^{(i)})$. 
The goal is therefore to construct an ordering of these velocity vectors that preserves proximity in velocity space. 

We construct this order using a Morton curve, a space-filling curve that interleaves the binary digits of integer grid coordinates to produce a single one-dimensional index \cite{morton_computer_1966, bader_space-filling_2013}. 
For a three-dimensional grid with $2^p$ grid points along each axis, the Morton index maps the integer triplet $(i_x, i_y, i_z)$ to an integer in $[0, 2^{3p}-1]$, approximately preserving spatial locality. 

The parameter $p$ is chosen such that the number of grid points exceeds the number of particles in the cell,
\begin{equation}
p = \min \{\,p \in \mathbb{N} : 2^{3p} > N_j\,\}.
\end{equation}
To map continuous velocity values to this discrete grid, each velocity component is first normalized to the unit interval,
\begin{equation}
M^{(i)}_{t,k,\mathrm{scaled}} 
= \frac{1}{2}\!\left[\,1 + \frac{M^{(i)}_{t,k} - \bar c_{t,k}}{3\sigma_{t,k}}\,\right],
\qquad k \in \{x, y, z\},
\end{equation}
where $\bar c_{t,k}$ and $\sigma_{t,k}$ denote the local mean velocity and standard deviation in the cell. 
Out-of-range values are clipped to the nearest bound.
The scaled components are then multiplied by $2^p$ and truncated to integers,
\begin{equation}
i_k^{(i)} = \operatorname{int}\!\left( 2^p\, M^{(i)}_{t,k,\mathrm{scaled}} \right),
\end{equation}
yielding the integer triplet $(i_x^{(i)}, i_y^{(i)}, i_z^{(i)})$ from which the Morton index is computed \cite{baert_libmorton_2018}. 

Sorting the Morton indices produces the desired one-dimensional ordering, in which nearby entries correspond to particles with similar velocities. 
To minimize computational overhead, radix sort is employed. 
Unlike comparison-based methods that require $O(N \log N)$ operations, radix sort proceeds by counting occurrences of digit values in a series of counting passes proportional to the number of digits $k$, resulting in an overall complexity of $O(kN)$ \cite{cormen_introduction_2022}.

Once the particle order is established, a corresponding set of quasi-random numbers is generated by producing $N_j$ three-dimensional points of a Sobol' sequence using Intel oneAPI MKL. 
To ensure unbiased estimators and enable variance estimation across independent runs while still maintaining low-discrepancy properties, the sequence is randomized by a digital shift operation, that is, a bitwise exclusive-or of each Sobol' point with a uniformly drawn binary vector \cite{dick_digital_2010}. 

The uniform quasi-random samples are then sorted according to the Morton order established from the velocities. 
Sorting the quasi-random samples rather than the particles achieves the same effect without modifying particle data.

Each ordered quasi-random number $u^{(i)} \in [0,1]$ is transformed into a standard normal variable by inverse transform sampling,
\begin{equation}
\xi^{(i)} = \Phi^{-1}(u^{(i)}),
\end{equation}
where $\Phi^{-1}$ denotes an approximation of the inverse cumulative distribution function of the standard normal distribution. 
The resulting $\xi^{(i)}$ replace the pseudo-random numbers in the Ornstein-Uhlenbeck velocity update discussed previously.

\subsection{Boundary Conditions}
Here, we focus on diffusive boundary conditions. 
When a particle collides with a wall, its velocity is resampled from a half-range Maxwellian distribution defined by the wall temperature and mean velocity \cite{cercignani_boltzmann_1988, bird_molecular_1994}. 

A uniform random number $U \sim \mathcal{U}(0,1)$ is used to generate the component of the velocity normal to the wall,
\begin{equation}
M_\text{normal} = \sqrt{-2\log(U)},
\end{equation}
and a Gaussian random variable $\xi \sim \mathcal{N}(0,1)$ provides the tangential component,
\begin{equation}
M_\text{tangent} = \xi.
\end{equation}
Since quasi-random numbers preserve their low-discrepancy properties only when generated in sequence, the wall reflections use pre-generated sequences assigned to each boundary cell.

\section{Numerical Results}
This section presents numerical experiments to assess how the FP--Array-RQMC method performs in comparison with alternative sampling strategies, in terms of error and convergence of estimated macroscopic quantities. We first introduce the sampling techniques considered for comparison, followed by the error metric used throughout this section. Four test cases are then examined. The first is a homogeneous relaxation toward a prescribed reservoir temperature, in which the drift and diffusion coefficients are fixed. The second is a homogeneous McKean--Vlasov-type relaxation, where these coefficients depend on ensemble moments. The third test case is an inhomogeneous Couette flow, which evaluates the method in a multi-cell setting with boundary interactions. The fourth is an inhomogeneous heat-flux configuration, included as an additional inhomogeneous benchmark for error comparison.

\subsection{Sampling Strategies}
We consider several sampling strategies that are commonly used in practice or that are natural candidates for variance reduction in the FP method. As a baseline, we employ standard pseudo-random sampling. A simple and commonly used improvement is the normalization of the random variables at each time step, such that their empirical mean is zero and their variance is one.

Another commonly used variance-reduction technique for stochastic differential equations is antithetical sampling. Random variables are generated in pairs $(\xi,-\xi)$, which enforces symmetry of the sampling distribution. In homogeneous relaxation problems, entire particle trajectories can be constructed antithetically, since both trajectories remain subject to identical drift and diffusion coefficients. In inhomogeneous configurations, particle transport between cells and boundary interactions break this symmetry, so that only the benefit of a symmetric sampling distribution remains.

In \cite{gorji_variance_2015}, a variance-reduction strategy based on a control variate formulation for the FP method was investigated. Let $M$ denote the particle velocity in the target process and $M_c$ the corresponding velocity in a control process with analytically known moments. For a moment of interest $\mu=\mathbb{E}[g(M)]$, we use the estimator
\begin{equation}
\hat{\mu}_{CV}
= \frac{1}{N}\sum_{i=1}^N \Big(g(M^{(i)})-g(M_c^{(i)})\Big) + \mu_c,
\qquad \mu_c=\mathbb{E}[g(M_c)] .
\end{equation}
The two processes are constructed using identical stochastic input, such that $g(M)$ and $g(M_c)$ are strongly correlated. As a result, the variance of the difference $g(M)-g(M_c)$ is significantly reduced.

Since the control process can be formulated most easily for a homogeneous test case with fixed drift and diffusion coefficients, we restrict the control variate strategy to this configuration.

Besides the Array-RQMC formulation, we also use quasi-random numbers that are shuffled at each time step. This preserves the high-quality quasi-random distribution while breaking correlations between successive time steps.

\subsection{Error Metrics}
To compare the different sampling strategies, we quantify the error of estimated macroscopic quantities, expressed in dimensionless form, with respect to a reference solution. Let $\mu^{(r)}_{t,j}$ denote the value of a macroscopic quantity in repetition $r$, at time step $t$ and spatial cell $j$, and let $\mu_{\mathrm{ref},t,j}$ denote the corresponding reference value.

For each $(t,j)$, we first compute the ensemble mean
\begin{equation}
\bar\mu_{t,j} = \frac{1}{R}\sum_{r=1}^{R} \mu^{(r)}_{t,j},
\end{equation}
and the ensemble variance
\begin{equation}
\mathrm{Var}_{t,j} = \frac{1}{R}\sum_{r=1}^{R} \big(\mu^{(r)}_{t,j}\big)^2 - \bar\mu_{t,j}^{\,2}.
\end{equation}

The bias with respect to the reference solution is defined as
\begin{equation}
\mathrm{Bias}_{t,j} = \big|\bar\mu_{t,j} - \mu_{\mathrm{ref},t,j}\big| .
\end{equation}

The root-mean-square error is then obtained from the standard bias--variance decomposition,
\begin{equation}
\mathrm{RMSE}_{t,j}
= \sqrt{\mathrm{Bias}_{t,j}^{\,2} + \mathrm{Var}_{t,j}} .
\end{equation}

To obtain a single scalar error measure for each simulation configuration, the RMSE is averaged over all time steps and spatial cells,
\begin{equation}
\overline{\mathrm{RMSE}}
= \frac{1}{N_{\mathrm{cells}}\,T_{\mathrm{steps}}}
  \sum_{t,j} \mathrm{RMSE}_{t,j}.
\end{equation}

This time--space averaged RMSE is used as the error measure in all convergence studies.

\subsection{Homogeneous Relaxation with Constant Coefficients}
We consider a spatially homogeneous, weakly rarefied monoatomic argon gas, initialized at equilibrium with temperature $T_0=300\,\mathrm{K}$ and number density $n=10^{19}\,\mathrm{m^{-3}}$, corresponding to a Knudsen number of $\mathrm{Kn}=0.17$. The system relaxes toward a reservoir temperature $T_\infty=600\,\mathrm{K}$ under constant drift and diffusion coefficients chosen accordingly.

The relaxation is simulated for $T_{\mathrm{steps}} = 35$ time steps with a time step size of $\Delta t = 10^{-4}$, which is approximately one tenth of the relaxation time $\tau$. By this time, the system is close to equilibrium. Each configuration is repeated $R = 1000$ times. The number of particles is chosen as $N = 2^{p}$ for $p = 6, \dots, 20$.

The sampling strategies considered are standard pseudo-random sampling, normalized pseudo-random sampling, antithetical sampling, a control variate formulation, quasi-random sampling with shuffling, and quasi-random sampling with the Array-RQMC technique.

The analytical solution of this relaxation serves as the reference. The mean velocity and the off-diagonal stress components remain constant in time, while the translational energy relaxes exponentially toward equilibrium with rate $1/\tau$,
\begin{equation}
c_k(t)=c_k(0), \qquad 
\sigma_{ij}(t)=\sigma_{ij}(0)\;\;(i\neq j), \qquad 
\epsilon(t)=\epsilon_\infty + \big(\epsilon(0)-\epsilon_\infty\big)e^{-t/\tau}.
\end{equation}

\begin{figure}[h]
\centering
\includegraphics[width=1\linewidth]{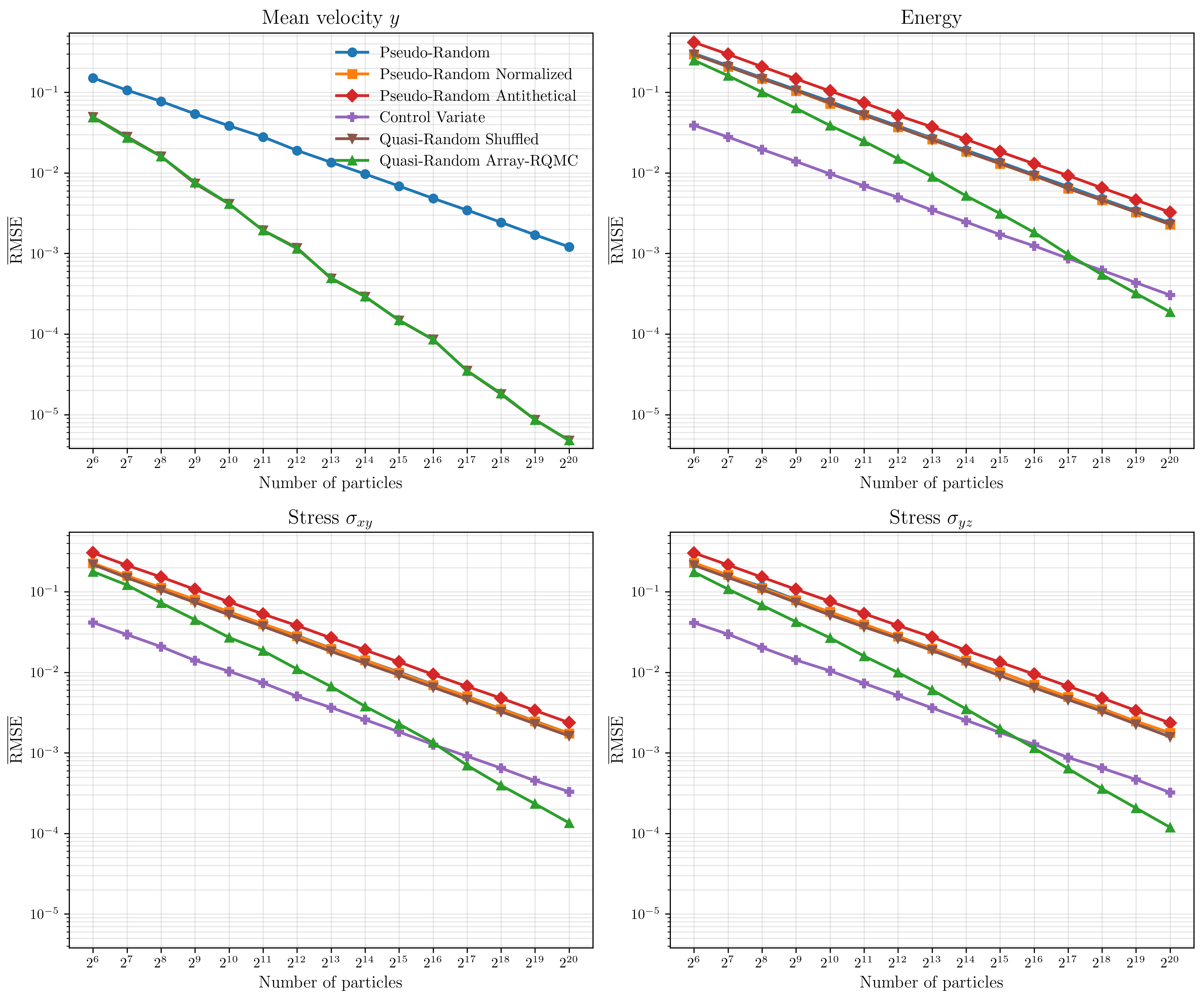}
\caption{Convergence of the average RMSE for a relaxation with constant coefficients. The four panels show the errors for the mean velocity component $y$, the translational energy, and the $\sigma_{xy}$ and $\sigma_{yz}$ components of the stress tensor. The figure compares pseudo-random sampling, normalized pseudo-random sampling, antithetical sampling, control variate, shuffled quasi-random sampling, and quasi-random sampling with the Array-RQMC technique.}
\label{fig:Constant_Coefficients_Relaxation_Convergence}
\end{figure}

Figure~\ref{fig:Constant_Coefficients_Relaxation_Convergence} shows the convergence of the average RMSE for several macroscopic quantities. All pseudo-random--based sampling strategies exhibit the expected convergence rate of $-0.5$, but differ significantly in their error magnitudes.

Normalized sampling, antithetic sampling, and the control variate formulation preserve a constant mean velocity -- normalized and antithetic sampling through symmetry of the sampling distribution, and the control variate formulation through cancellation between the target and control process contributions in the estimator. Since the initial distribution is normalized to zero mean, this constant mean results in a vanishing error for the mean velocity.

For the second-order moments, normalization does not reduce the error. Antithetical sampling yields errors larger by approximately a factor of two, implying that twice as many particles are required to achieve the same accuracy as with standard pseudo-random sampling. The variance-reduction mechanism of antithetical sampling relies on a negative correlation between paired sample contributions \(g(M)\) and \(g(-M)\) inside the estimator. However, for quadratic observables, \(g(M)=g(-M)\), leading to a perfect positive correlation. Consequently, the estimated moments depend only on the independent half of the samples, which explains the observed factor-of-two increase in error.

\begin{table}[h!]
\centering
\caption{Convergence rates of the average RMSE for the homogeneous relaxation with constant coefficients.}
\label{tab:Constant_Coefficients_Relaxation_Convergence}
\begin{tabular}{lcccc}
\toprule
Sampling strategy 
& Mean velocity $y$ 
& Energy 
& Stress $\sigma_{xy}$ 
& Stress $\sigma_{yz}$ \\
\midrule
Pseudo-random                   & $-0.50$ & $-0.50$ & $-0.50$ & $-0.50$ \\
Pseudo-random Normalized        & $0$     & $-0.50$ & $-0.50$ & $-0.50$ \\
Pseudo-random Antithetical      & $0$     & $-0.50$ & $-0.50$ & $-0.50$ \\
Pseudo-random Control variate                 & $0$     & $-0.50$ & $-0.50$ & $-0.50$ \\
Quasi-random Shuffled            & $-0.96$ & $-0.50$ & $-0.50$ & $-0.50$ \\
Quasi-random Array-RQMC         & $-0.96$ & $-0.75$ & $-0.75$ & $-0.75$ \\
\bottomrule
\end{tabular}
\end{table}

In contrast, the control variate formulation yields a substantial reduction of the error for the second-order macroscopic quantities by exploiting a strongly correlated control process with analytically known moments. However, the applicability of this approach depends on the availability of a suitable control process.

Both quasi-random sampling strategies yield a significantly smaller error and a steeper convergence rate of approximately $-0.96$ for the mean velocity. Both methods perform identically, since the mean velocity is a linear observable and therefore the estimator depends only on the sample values and not on their ordering.

For the second-order moments, shuffled quasi-random sampling does not provide any advantage over standard pseudo-random sampling, because the low-discrepancy structure is lost across time steps.

For small particle numbers, the Array-RQMC method yields only a marginal improvement over pseudo-random sampling and performs worse than the control variate formulation. However, by preserving the low-discrepancy structure across time steps, Array-RQMC achieves a significantly improved convergence rate of approximately $-0.75$, so that the error becomes smaller than that of the control variate formulation when sufficiently many particles are used.

Table~\ref{tab:Constant_Coefficients_Relaxation_Convergence} summarizes the convergence rates for all quantities and sampling strategies.

\subsection{Homogeneous McKean--Vlasov Relaxation}
This relaxation test case uses drift and diffusion coefficients that depend on the velocity moments, corresponding to the standard McKean--Vlasov-type formulation of the FP equation.

To initialize the particles in a non-equilibrium state with a non-zero stress tensor, we design an anisotropic distribution by modifying a Maxwellian through a planar cut. The cutting plane passes through the origin and is oriented $25^\circ$ counterclockwise with respect to the $x$--$y$ plane. All velocities lying on the side of the plane corresponding to increasing $z$ are discarded. The remaining velocities are scaled to enforce zero mean and unit variance.

The gas relaxes toward a weakly rarefied equilibrium state with $T=300\,\mathrm{K}$ and $n=10^{19}\,\mathrm{m^{-3}}$, corresponding to a Knudsen number of $\mathrm{Kn}=0.17$.

The relaxation is simulated for $T_{\mathrm{steps}} = 100$ time steps with a time step size of $\Delta t = 10^{-4}$, corresponding to approximately one tenth of the relaxation time $\tau$. At this point, the system is close to equilibrium. Each configuration is repeated $R = 1000$ times. The number of particles is chosen as $N = 2^{p}$ for $p = 6, \dots, 11$. The sampling strategies considered are standard pseudo-random sampling, normalized pseudo-random sampling, antithetical sampling, quasi-random sampling with shuffling, and quasi-random sampling with the Array-RQMC technique.

The analytical solution of the relaxation is used as the reference. The mean velocity and the translational energy remain constant in time, while the components of the stress tensor relax exponentially with rate $1/\tau$,
\begin{equation}
c_k(t)=c_k(0), \qquad \epsilon(t)=\epsilon(0), \qquad \sigma_{ij}(t)=\sigma_{ij}(0)\,e^{-t/\tau}.
\end{equation}

\begin{figure}[h!]
\centering
\includegraphics[width=1\linewidth]{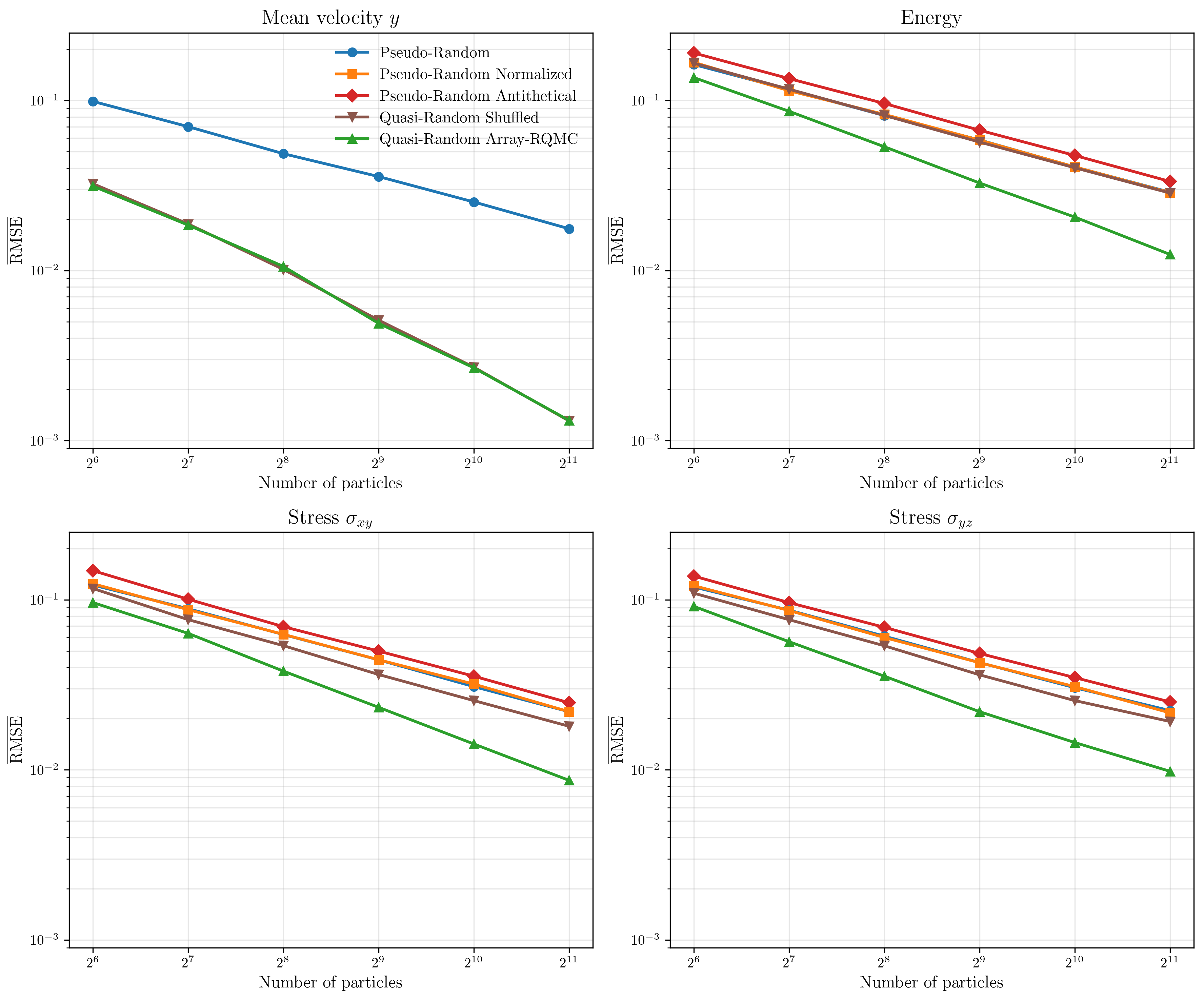}
\caption{Convergence of the average RMSE for the homogeneous McKean--Vlasov relaxation test case. The four panels show the errors for the mean velocity component $y$, the translational energy, and the $\sigma_{xy}$ and $\sigma_{yz}$ components of the stress tensor. The figure compares pseudo-random sampling, normalized pseudo-random sampling, antithetical sampling, shuffled quasi-random sampling, and quasi-random sampling with the Array-RQMC technique.}
\label{fig:relaxation_convergence}
\end{figure}

Figure~\ref{fig:relaxation_convergence} shows the convergence of the average RMSE for four macroscopic quantities. All pseudo-random--based sampling strategies exhibit the expected convergence rate of $-0.5$. Normalizing the pseudo-random samples reduces the error by a constant factor. Antithetical sampling yields slightly larger errors. Both normalized and antithetical sampling result in a vanishing error for the mean velocity.

\begin{table}[h!]
\centering
\caption{Convergence rates of the average RMSE for the homogeneous McKean--Vlasov relaxation test case.}
\label{tab:relaxation_slopes}
\begin{tabular}{lcccc}
\toprule
Sampling strategy 
& Mean velocity $y$ 
& Energy 
& Stress $\sigma_{xy}$ 
& Stress $\sigma_{yz}$ \\
\midrule
Pseudo-random                   & $-0.49$ & $-0.50$ & $-0.50$ & $-0.49$ \\
Pseudo-random Normalized        & $0$     & $-0.51$ & $-0.50$ & $-0.50$ \\
Pseudo-random Antithetical      & $0$     & $-0.50$ & $-0.51$ & $-0.49$ \\
Quasi-random Shuffled           & $-0.93$ & $-0.51$ & $-0.51$ & $-0.51$ \\
Quasi-random Array-RQMC         & $-0.92$ & $-0.69$ & $-0.70$ & $-0.65$ \\
\bottomrule
\end{tabular}
\end{table}

Table~\ref{tab:relaxation_slopes} summarizes the convergence rates for all quantities and sampling strategies.

Overall, the Array-RQMC formulation retains most of its advantages in the McKean--Vlasov relaxation compared to the relaxation with constant coefficients. Although the convergence rates are slightly reduced, the method continues to provide a substantial error reduction, demonstrating its robustness with respect to moment-dependent drift and diffusion.

\subsection{Inhomogeneous Couette Flow}
\begin{figure}[h]
\centering
\begin{tikzpicture}[scale=1.0, >=Stealth]

\draw[thick] (0,0) rectangle (6,2);
\foreach \x in {0.5,1.0,...,5.5} {
  \draw[densely dotted] (\x,0) -- (\x,2);
}
\draw[->] (-4,0.6) -- (-3,0.6) node[right] {$y$};
\draw[->] (-4,0.6) -- (-4,1.6) node[above] {$x$};
\node[left] at (-0.2,1.4) {$T_\text{left}=300\,\mathrm{K}$};
\node[left] at (0,0.8) {$v_\text{left}=0\,\mathrm{ms^{-1}}$};
\node at (3,1.2) {$T_0=300\,\mathrm{K}$};
\node at (3,0.7) {$v_0=0$};
\node[right, align=left] at (6,1.5) {$T_\text{right}=300/400\,\mathrm{K}$};
\node[right, align=left] at (6,0.9) {$v_\text{right}=100/0\,\mathrm{m\,s^{-1}}$};

\end{tikzpicture}
\caption{Schematic of the inhomogeneous test configurations. A monoatomic argon gas is initially at rest with $T_0=300\,\mathrm{K}$ between two parallel plates. The left plate is fixed at $T_\text{left}=300\,\mathrm{K}$ and $v_\text{left}=0$. The Couette-flow is obtained by moving the right plate with $v_\text{right}=100\,\mathrm{m\,s^{-1}}$ at constant temperature, while the heat-flux is obtained by heating the right plate to $T_\text{right}=400\,\mathrm{K}$ while keeping it stationary. Dotted lines indicate the spatial cell discretization.}
\label{fig:inhomogeneous_schematic}
\end{figure}
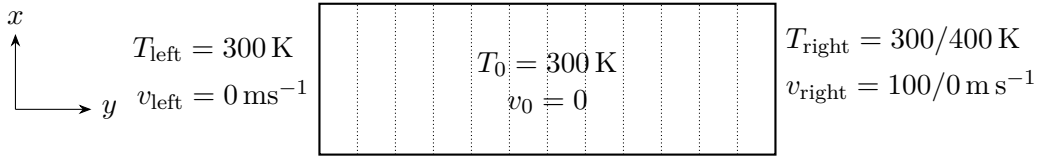

We consider the buildup of a Couette flow to evaluate the method in an inhomogeneous setting with multiple cells and boundary interactions (see Fig.~\ref{fig:inhomogeneous_schematic}). A weakly rarefied monoatomic argon gas is initially in equilibrium between two parallel plates at temperature $T=300\,\mathrm{K}$ and number density $n=10^{19}\,\mathrm{m^{-3}}$, corresponding to a Knudsen number of $\mathrm{Kn}=0.17$. At time $t=0$, the upper plate is set into motion with a velocity of $100\,\mathrm{m\,s^{-1}}$ in $y$-direction, while the lower plate remains at rest.

The buildup of the Couette flow is simulated with a time step size of $\Delta t = 10^{-4}$ for $T_{\mathrm{steps}} = 300$ time steps, such that a steady state is reached. Each configuration is repeated $R = 50$ times. The number of particles is chosen as $N = 2^{p}$ with $p = 6, \dots, 14$. The domain is discretized into 10 cells in the wall-normal direction. The sampling strategies considered are standard pseudo-random sampling, normalized pseudo-random sampling, antithetical sampling, quasi-random sampling with shuffling, and quasi-random sampling with the Array-RQMC technique.

Since no analytical solution is available for the buildup of a Couette flow of a rarefied gas, a high-resolution FP simulation with pseudo-random sampling is used as reference. The reference solution is computed with $N_{\mathrm{ref}} = 10^{6}$ particles and $R_{\mathrm{ref}} = 100$ independent repetitions. The macroscopic quantities are obtained by averaging over all repetitions.

\begin{figure}[h!]
\centering
\includegraphics[width=1\linewidth]{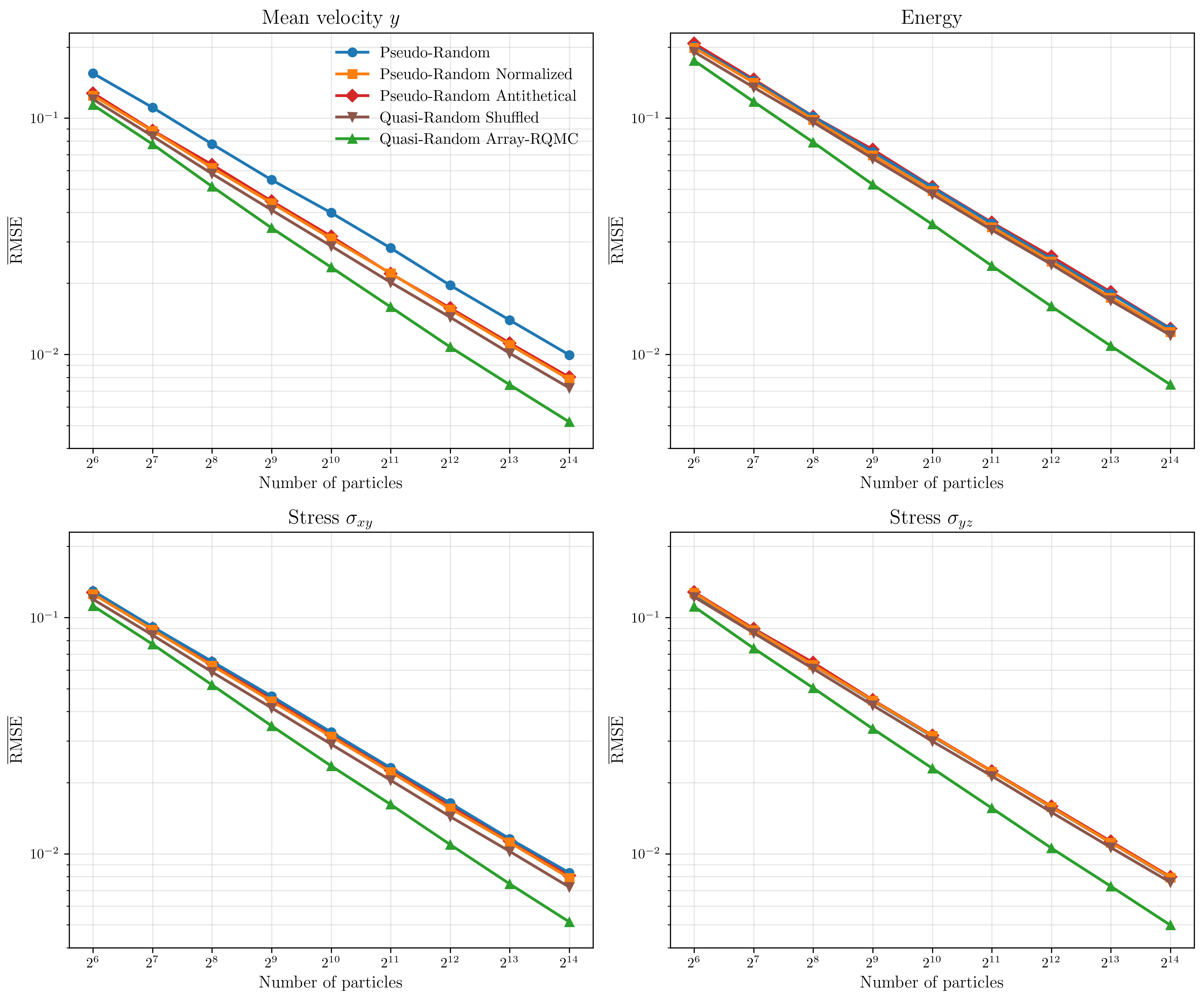}
\caption{Convergence of the average RMSE for the inhomogeneous Couette flow test case. The four panels show the errors for the mean velocity component $y$, the translational energy, and the $\sigma_{xy}$ and $\sigma_{yz}$ components of the stress tensor. The figure compares pseudo-random sampling, normalized pseudo-random sampling, antithetical sampling, shuffled quasi-random sampling, and quasi-random sampling with the Array-RQMC technique.}
\label{fig:couette_convergence}
\end{figure}

Figure~\ref{fig:couette_convergence} shows the convergence of the average RMSE for the inhomogeneous Couette flow. All pseudo-random--based sampling strategies exhibit the expected convergence rate of approximately $-0.5$. Normalized and antithetical sampling yield a comparable constant-factor reduction of the error for the mean velocity. In contrast to the homogeneous relaxation cases, a symmetric sampling distribution no longer leads to a vanishing error in the mean, since the mean velocity is continuously perturbed by particle transport between cells and by boundary interactions.

For the second-order moments, all pseudo-random sampling strategies, including antithetical sampling, produce nearly identical errors. This contrasts with the homogeneous setting, where antithetical sampling resulted in larger errors. Although the additional randomness introduced by particle transport and boundary interactions breaks up the antithetical structure over time, it also compensates for the reduced number of independent samples.

\begin{table}[h!]
\centering
\caption{Convergence rates of the average RMSE for the inhomogeneous Couette flow test case.}
\label{tab:couette_slopes}
\begin{tabular}{lcccc}
\toprule
Sampling strategy 
& Mean velocity $y$ 
& Energy 
& Stress $\sigma_{xy}$ 
& Stress $\sigma_{yz}$ \\
\midrule
Pseudo-random                   & $-0.50$ & $-0.50$ & $-0.50$ & $-0.50$ \\
Pseudo-random Normalized        & $-0.50$ & $-0.50$ & $-0.50$ & $-0.50$ \\
Pseudo-random Antithetical      & $-0.50$ & $-0.50$ & $-0.50$ & $-0.50$ \\
Quasi-random Shuffled           & $-0.51$ & $-0.50$ & $-0.51$ & $-0.50$ \\
Quasi-random Array-RQMC         & $-0.56$ & $-0.57$ & $-0.56$ & $-0.56$ \\
\bottomrule
\end{tabular}
\end{table}

In contrast to the homogeneous relaxation cases, shuffling the quasi-random numbers does not improve the convergence rate for the mean velocity and yields errors comparable to normalized and antithetical pseudo-random sampling. In the homogeneous case, the linearity of the estimator implies that the mean velocity depends only on the quasi-random point set, while the ordering of the samples is irrelevant. In the inhomogeneous case, however, the ordering of samples becomes important, because particle transport and boundary interactions introduce additional sources of randomness that couple successive updates. Although quasi-random point sets are also used for particle initialization and boundary sampling, the interplay of these point sets breaks the global low-discrepancy structure when the samples are shuffled. For the second-order moments, shuffled quasi-random numbers show no improvement over pseudo-random sampling.

The Array-RQMC formulation retains an advantage in the inhomogeneous setting. Although the additional randomness reduces its effectiveness compared to the homogeneous case, Array-RQMC still preserves part of the low-discrepancy structure across time steps. As a result, it achieves an improved convergence rate of approximately $-0.56$ for all macroscopic quantities, as summarized in Table~\ref{tab:couette_slopes}. The error reduction is marginal for small particle numbers, but becomes increasingly significant as $N$ grows. Already for $N=2^{13}$ particles, Array-RQMC is approximately twice as efficient as pseudo-random sampling.

Overall, the additional randomness inherent in the inhomogeneous settings drives all sampling strategies closer to the pseudo-random baseline. Antithetical sampling benefits from this behavior, while the performance of Array-RQMC is weakened. Nevertheless, Array-RQMC preserves a clear advantage, which becomes increasingly significant as the number of particles grows. This demonstrates that low-discrepancy structure can still be exploited in inhomogeneous problems, even though its impact is diminished by transport and boundary interactions.

\section{Inhomogeneous Heat Flux}
We consider the buildup of a heat flux to further evaluate the method in an inhomogeneous setting with spatial transport and boundary interactions (see Fig.~\ref{fig:inhomogeneous_schematic}). A weakly rarefied monoatomic argon gas is initially in equilibrium between two parallel plates at temperature $T=300\,\mathrm{K}$ and number density $n=10^{19}\,\mathrm{m^{-3}}$, corresponding to a Knudsen number of $\mathrm{Kn}=0.17$. At time $t=0$, the temperature of the upper plate is raised to $400\,\mathrm{K}$, while the lower plate remains at $300\,\mathrm{K}$, leading to the development of a heat flux in the wall-normal direction.

The heat-flux buildup is simulated with a time step size of $\Delta t = 2\times10^{-4}$ for $T_{\mathrm{steps}} = 150$ time steps, until a steady state is reached. Each configuration is repeated $R = 50$ times. The number of particles is chosen as $N = 2^{p}$ with $p = 6, \dots, 17$. The domain is discretized into 20 cells in the wall-normal direction. The sampling strategies considered are standard pseudo-random sampling, normalized pseudo-random sampling, antithetical sampling, quasi-random sampling with shuffling, and quasi-random sampling with the Array-RQMC technique.

Since no analytical solution is available for the buildup of a heat flux in a rarefied gas, a high-resolution FP simulation with pseudo-random sampling is used as reference. The reference solution is computed with $N_{\mathrm{ref}} = 10^{6}$ particles and $R_{\mathrm{ref}} = 100$ independent repetitions. The macroscopic quantities are obtained by averaging over all repetitions.

\begin{figure}[h!]
\centering
\includegraphics[width=1\linewidth]{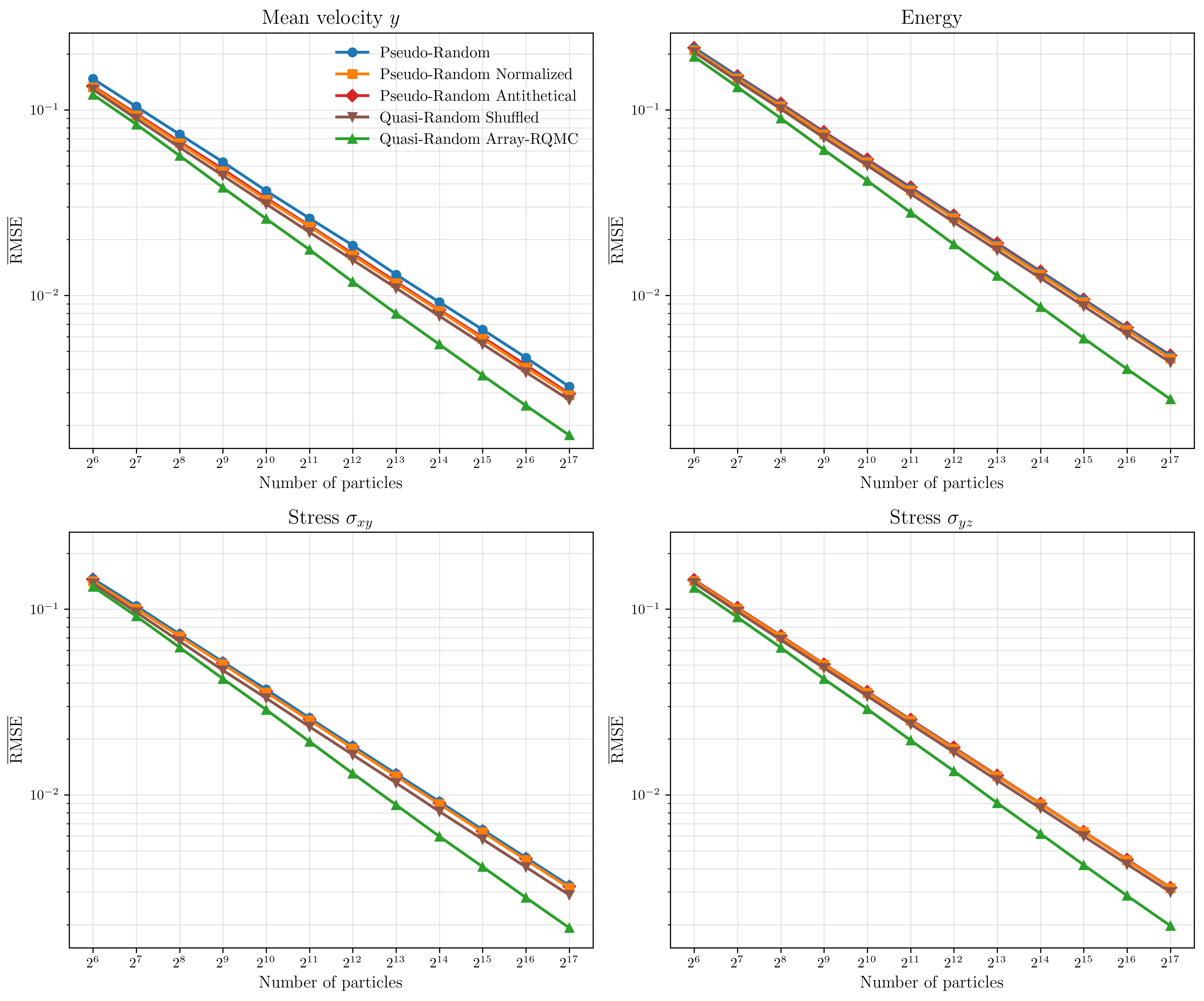}
\caption{Convergence of the average RMSE for the inhomogeneous heat flux test case. The four panels show the errors for the mean velocity component $y$, the translational energy, and the $\sigma_{xy}$ and $\sigma_{yz}$ components of the stress tensor. The figure compares pseudo-random sampling, normalized pseudo-random sampling, antithetical sampling, shuffled quasi-random sampling, and quasi-random sampling with the Array-RQMC technique.}
\label{fig:heat_flux_convergence}
\end{figure}

Figure~\ref{fig:heat_flux_convergence} shows the convergence of the average RMSE for the inhomogeneous heat-flux configuration. All pseudo-random--based sampling strategies exhibit the expected convergence rate of approximately $-0.5$. Normalized and antithetical sampling yield a marginal constant-factor reduction of the error for the mean velocity, which is even weaker than in the Couette-flow case. For the second-order moments, all pseudo-random sampling strategies produce essentially identical errors.

Shuffled quasi-random sampling does not exhibit any improvement over pseudo-random sampling. Its convergence rates and error magnitudes are essentially identical to those of the pseudo-random--based strategies.

The Array-RQMC formulation again retains a moderate advantage in the heat-flux configuration. Despite the increased spatial resolution and the larger time step, the convergence rate remains close to $-0.56$ for all macroscopic quantities, as summarized in Table~\ref{tab:heat_flux_slopes}. As in the Couette-flow case, the error reduction is marginal for small particle numbers but becomes increasingly significant as $N$ grows. Already for $N=2^{14}$ particles, Array-RQMC is approximately twice as efficient as pseudo-random sampling.

\begin{table}[h!]
\centering
\caption{Convergence rates of the average RMSE for the inhomogeneous heat-flux test case.}
\label{tab:heat_flux_slopes}
\begin{tabular}{lcccc}
\toprule
Sampling strategy 
& Mean velocity $y$ 
& Energy 
& Stress $\sigma_{xy}$ 
& Stress $\sigma_{yz}$ \\
\midrule
Pseudo-random                   & $-0.50$ & $-0.50$ & $-0.50$ & $-0.50$ \\
Pseudo-random Normalized        & $-0.50$ & $-0.50$ & $-0.50$ & $-0.50$ \\
Pseudo-random Antithetical      & $-0.50$ & $-0.50$ & $-0.50$ & $-0.50$ \\
Quasi-random Shuffled           & $-0.50$ & $-0.50$ & $-0.51$ & $-0.50$ \\
Quasi-random Array-RQMC         & $-0.56$ & $-0.56$ & $-0.56$ & $-0.55$ \\
\bottomrule
\end{tabular}
\end{table}

\section{Conclusion}
We integrated quasi-random numbers into the FP method using the Array-RQMC technique, which preserves the low-discrepancy structure across time steps by sorting the particles in velocity space. Since this mechanism is local in both time and space, the FP--Array-RQMC formulation can exploit the advantages of quasi-random numbers in spatially discretized settings with particle transport across cells and boundary interactions.

Numerical experiments in both homogeneous and inhomogeneous test cases demonstrate that FP--Array-RQMC improves the convergence behavior of macroscopic estimators. In homogeneous relaxation problems, convergence rates increase to approximately $-0.9$ for first-order moments and to about $-0.7$ for second-order moments. In inhomogeneous problems, the coupling of velocity updates, particle transport across cells, and boundary interactions introduces additional randomness, despite the low-discrepancy structure of each individually. This interplay reduces the convergence rates, but a moderate improvement to about $-0.56$ is still preserved. Because techniques based on pseudo-random sampling do not change the convergence rate, the benefit of FP--Array-RQMC increases with particle number, resulting in substantially smaller estimator errors for large $N$.

Future work will aim to improve the coupling between velocity updates, particle transport, and boundary interactions in order to further enhance convergence rates in inhomogeneous problems, for example through kernel-based strategies. In addition, the FP--Array-RQMC approach will be explored for more advanced FP models matching relaxation rates of higher-order moments, as well as polyatomic gases.

\newpage
\section*{Acknowledgements}
This work was supported by the Deutsche Forschungsgemeinschaft (DFG) within the Priority Program SPP 2410 ``Hyperbolic Balance Laws in Fluid Mechanics: Complexity, Scales and Randomness (CoScaRa)'', under project number 525660607 ``Model Cascades for Stochastic Particle Simulations of Rarefied Polyatomic Gases''.

\bibliographystyle{plain}
\bibliography{library}

\end{document}